\documentclass[11pt,leqno,a4paper]{article}

\usepackage{amsfonts}
\usepackage{amsmath,amssymb,amsthm}

\theoremstyle{plain}
\newtheorem{thm}{Theorem}
\newtheorem{lem}{Lemma}[section]
\newtheorem{cor}{Corollary}[section]
\renewcommand{\d}{\displaystyle}

\DeclareMathOperator*{\Res}{Res}

\numberwithin{equation}{section}

\begin{document}
\title{Conditional Results for a Class of Arithmetic Functions:
a variant of  H. L. Montgomery and R. C. Vaughan's method }
\author{Xiaodong Cao and Wenguang Zhai }
\date{}

\footnotetext[0]{2000 Mathematics Subject Classification: 11N37,
11L07.}

\footnotetext[0]{Key Words:  square-full integer, exponential
convolution, exponential divisor, e-$r$-free integer,
$(r,l)$-integer, exponent pair. }

\footnotetext[0]{This work is supported by the National Natural
Science Foundation of China(Grant No. 11171344) and the Natural
Science Foundation of Beijing(Grant No. 1112010).}
 \maketitle
 
{\bf Abstract.} Let $a, b ,c $ and $k$ be positive integers such that $1\leq a\leq
b,a<c<2(a+b), c\ne b$ and $(a,b,c)=1$. Define the  arithmetic function $f_k(a,b;c;n)$  by
$$
\sum_{n=1}^{\infty}\frac{f_k(a,b;c;n)}{n^s}=\frac{\zeta (as)\zeta
(bs)}{\zeta^k(cs)}, \Re s >1.$$
 Let $\Delta_k(a,b;c;x)$ denote the error term of the summatory 
 function of the function $f_k(a,b;c;n).$ IN this paper we shall give two 
 expressions of  $\Delta_k(a,b;c;x)$. As applications, we study the so-called 
 $(l,r)$-integers, the generalized square-full integers, the $e-r$-free integers,
 the divisor problem over $r$-free integers, the $e$-square-free integers. 
 An important tool is a  generalization of a method of H. L. Montgomery and R. C. Vaughan.

\section{Introduction and main results}
   W. G. Nowak\cite{nowak}, M. K\"{u}leitner and W. G. Nowak \cite{kn} studied a class of very general arithmetic function $a(n),$ which possess a generating Dirichlet series
 $$
\sum_{n=1}^{\infty}\frac {a(n)}{n^s}=\frac{f_1(m_1s)\cdots f_K(m_Ks)}{g_1(n_1s)\cdots g_J(n_Js)} h(s),
 $$
 where $f_k$ and $g_j$ are certain generalizations of Riemann zeta-function,  $m_1\le \cdots \le m_K$ and $n_1\le \cdots \le n_J$
are natural numbers, and $h(s)$ is a good function which is regular and bounded in a sufficiently large half-plane. People are usually concerned with the summatory function $\sum_{n\leq x}a(n),$ especially sharp upper and lower bounds of its error term.  The above two papers give an upper bound and a lower bound for $a(n)$ in a very general sense.   Some special cases are also studied, see for example, \cite{ba1, ba2, ku, wu2}.

 The aim of this paper is to study a special case of $a(n),$ in which case we can get better upper results. Let $a, b ,c $ and $k$ be positive integers such that $1\leq a\leq
b,a<c<2(a+b), c\ne b$ and $(a,b,c)=1$. Let $\zeta(s)$ denote the
Riemann zeta-function. The arithmetic function $f_k(a,b;c;n)$ is
defined by
\begin{eqnarray}
\sum_{n=1}^{\infty}\frac{f_k(a,b;c;n)}{n^s}=\frac{\zeta (as)\zeta
(bs)}{\zeta^k(cs)}, \Re s >1.
\end{eqnarray}
In this paper we are concern with  the summatory function
\begin{eqnarray}
A_k(a,b;c;x):=\sum_{n\le x}f_k(a,b;c;n),\ x\geq 2.
\end{eqnarray}
The expected asymptotic formula of $A_k(a,b;c;x)$ is of the form
\begin{eqnarray}
 A_k(a,b;c;x) =\frac {\zeta (\frac ba)}{\zeta^k
(\frac ca)}x^{\frac 1a}+\frac {\zeta (\frac ab)}{\zeta^k (\frac
cb)}x^{\frac 1b}+\Delta_k(a,b;c;x)
\end{eqnarray}
when $a\ne b.$ When $a = b$, then an appropriate limit should be
taken in the above formula. As usual, $\Delta_k(a,b;c;x)$ is called
the error term of the function $A_k(a,b;c;x).$
 For convenience, we also use
notations $A(a,b;c;x)$, $\Delta(a,b;c;x)$ to denote
$A_1(a,b;c;x)$, $\Delta_1(a,b;c;x)$, respectively.

  By  Theorem 2 of M.
K\"{u}leitner and W. G. Nowak\cite{kn} or Theorem 3 of W. G. Nowak\cite{nowak}, it is easy to prove that
\begin{eqnarray}
\Delta_k(a,b;c;x):=\Omega\left(x^{\max (\frac
{1}{2(a+b)},\frac{1}{2c})}\right).
\end{eqnarray}
Hence one may conjecture that
\begin{eqnarray}
\Delta_k(a,b;c;x):=O\left(x^{\max (\frac
{1}{2(a+b)},\frac{1}{2c})+\varepsilon}\right).
\end{eqnarray}

    Many special cases of the function $f_k(a,b;c;n)$ have been
extensively studied in  number theory. We take some examples.

(1) The case $(a,b,c,k)=(2,3,6,1)$ is the well-known square-full
number problem(see \cite{ca,co2,cd,ns,wu2}). Suppose $a\nmid b,$
the cases $(a,b,c,k)=(a,b,2b,1)$ or $(a,b,c,k)=(b,a,2b,1)$ are
studied for the the generalized  square-full number problem(see
\cite{co1,su2} ).

(2) Suppose $r\geq 2$ is a fixed integer, the case
$(a,b,c,k)=(1,1,r,1)$ is the $r$-free divisor problem(see
\cite{ba1,ba2, fz1,fz2,ku}).

(3) Suppose $1<r<l$ are fixed integers, the case
$(a,b,c,k)=(1,r,l,1)$ corresponds to the distribution of the
so-called $(r,l)$-integers(see \cite{co4,sh,ss1,ss2,su3,ya}).

(4) Suppose $r\geq 1$ is a fixed integer, the case
$(a,b,c,k)=(1,2^r+1,2^r,1)$ corresponds to the the distribution of
the so-called e-$r$-free integers(see \cite{cz2,su1,to1,to2,wu1}).

(5) Suppose $r\geq 2$ is a fixed integer, the case
$(a,b,c,k)=(1,1,r,r+1)$ corresponds to the Dirichlet divisor
problem over the set of $r$-free integers(see \cite{fen,re}).

\bigskip

From the right-hand side of (1.1) it is  easily seen that the
unconditional asymptotic formula we could possibly prove at present
is at most
\begin{eqnarray}
 A_k(a,b;c;x) =\frac {\zeta (\frac ba)}{\zeta^k
(\frac ca)}x^{\frac 1a}+\frac {\zeta (\frac ab)}{\zeta^k (\frac
cb)}x^{\frac 1b}+O(x^{\frac 1c}\exp(-A (\log x)^{\frac 35}(\log
\log x)^{-\frac 15})),
\end{eqnarray}
where $A>0$ is some absolute constant.  Now define
$\theta_k(a,b;c)$ denote the infimum of $\alpha_k(a,b;c)$ such
that the estimate
\begin{equation}\Delta_k(a,b;c;x)\ll x^{\alpha_k(a,b;c)+\varepsilon}\end{equation}
holds.

 From (1.1) we also see that the function $f_k(a,b;c;n)$ is
related to the divisor function $d(a,b;n):=\sum_{n=m_1^am_2^b}1,$
which satisfies
$$\sum_{n=1}^\infty \frac{d(a,b;n)}{n^s}=\zeta(as)\zeta(bs), \Re
s>1.$$
We write
\begin{eqnarray}
D(a,b;x):= \sum_{n\leq x}d(a,b;n)=\zeta (\frac ba)x^{\frac
1a}+\zeta (\frac ab)x^{\frac 1b} +\Delta(a,b;x)
\end{eqnarray}
for $a\ne b$, and let
  $0<\alpha(a,b)<1/(a+b)$ be a
real number such that the estimate
\begin{equation}
\Delta(a,b;x)\ll x^{\alpha(a,b)+\varepsilon}
\end{equation} holds.

  As usual, $\Delta(a,b;x)$ is called the error
term of the asymmetric two-dimensional divisor problems.  For the history and classical results of $\Delta(a,b;x),$ see for
example\cite{iv,ikkn,kr,mn,ri,zc1}.

If $\alpha(a,b)\geq 1/c,$ then by the convolution approach we get
easily   $\theta_k(a,b;c)\leq  \alpha(a,b).$ Thus the difficulty of the  evaluation
of the function $A_k(a,b;c;x)$ is basically the difficulty of the
evaluation of the function $D(a,b;x).$ Without the loss of
generality, we always suppose later that $\alpha(a,b)< 1/c.$



The exponent $1/c$ in the error term in (1.6) is  closely related
to the distribution of the non-trivial zeros of $\zeta(s).$ People
usually assume the Riemann-hypothesis (RH) to reduce the constant
$1/c.$ See for example, \cite{ba1,ns,ss2,su2}. From now on, we
always suppose that RH holds.

In 1981, Montgomery and
Vaughan\cite{mv} developed a new ingenious method to treat the distribution of
$r$-free integers, which was also used by
many other authors, see for example, Baker\cite{ba1}, Nowak and
Schmeier\cite{ns}, Nowak\cite{nowak} etc.
However, as  W. G.
Nowak and M. Schmeier\cite{ns} observed
 in subsection (The Divisor Problem For
$(l,r)$-Integers) that: The $r=2, l=3$ is some exceptional.
 That is, in some cases, by Montgomery-Vaughan's method one could not
get directly better estimates than the usual  approach.

The main aim of this paper is to find a suitable expression of the
error term in (1.3) for every case. We have to consider two different cases: $c>b$ and
$c<b$. For these two cases, we have to use different convolution approaches.

Consider first   $c>b.$ We define the function $\mu_k$ by
\begin{eqnarray}
\sum_{n=1}^{\infty}\frac {\mu_k(n)}{n^s}=\frac {1}{\zeta ^k(s)},
\Re s>1.
\end{eqnarray}
(Also see Titchmarsh\cite{ti}, page 165-166.) Clearly $\mu_1$ is
the well-known M\"obius function $\mu$. Then Theorem 2 of Nowak\cite{nowak} essentially implies the following theorem.

\begin{thm}(W. G. Nowak) Let $x\ge 2,$ $a\le b<c<2(a+b)$ and $\Delta_k(a,b;c;x)$ be
defined by (1.3). If the RH is true, then for any $ 1\le y<
x^{\frac 1c}$
\begin{eqnarray}
\Delta_k(a,b;c;x)=\sum_{l\le y}\mu_k (l)\Delta\left(a,b;\frac
{x}{l^c}\right)+O\left(x^{\frac {1}{2a}+\varepsilon}y^{\frac
12-\frac {c}{2a}}+x^\varepsilon\right).
\end{eqnarray}
\end{thm}

Taking $y=x^{\frac {1-2a\alpha (a,b)}{a+c-2ac\alpha (a,b)}}$ and
noting  $\mu_{k} (l)\ll l^\varepsilon$,  we get

\begin{cor}   Suppose RH. If  $a\leq b<c<2(a+b)$ and
$\alpha{(a,b)}<\frac 1c$, then
\begin{eqnarray}
\theta_k(a,b;c)\leq  \frac {1-a\alpha (a,b)}{a+c-2ac\alpha (a,b)}.
\end{eqnarray}
\end{cor}

{\bf Remark 1.1.}  The
proof of Theorem 1 is based on a classical idea of Montgomery and
Vaughan\cite{mv} and the Dirichlet
convolution
\begin{equation}
f_k(a,b;c;n)=\sum_{n=m_1m_2^c}d(a,b;m_1) \mu_k(m_2).
\end{equation}

{\bf Remark 1.2.} Very fortunately, for arithmetical function $\mu_k $ we have an analogue of the well-known Vaughan's identity of M\"obius function $\mu$(see Lemma 4.1 below, in fact this is the third useful Vaughan-type's identity except the well-known von Manlgoldt function $\Lambda$ and M\"obius function $\mu$ ). By the method of exponential sums, one could
improve the result in Corollary 1.1.  For the related works, we refer
to papers \cite{ba1,ba2,ca,ku,wu2}. However, this is not the main
aim of the present paper.

\bigskip

We now turn to  the case $c<b.$ Some examples of this type can be found in \cite{ss1,ss2,su2,su3}. In this case, we hope to find an
estimate of the form  $\Delta_k(a,b;c;x)\ll x^{\alpha}\ (\alpha
<1/b)$ such that the second main term $\frac {\zeta (\frac
ab)}{\zeta^k (\frac cb)}x^{\frac 1b}$ becomes a real main term.  We
can also use the convolution (1.13) as our first choice to study
$A_k(a,b;c;x)$. Actually it is easy to check that Theorem 1 also
holds for $b/2<c<b.$ But when $a+c\le b,$ we have checked that it is very difficult to
prove $\theta_k(a,b;c)<1/b$ via Theorem 1 directly(Also see page 9, section 5 in \cite{cz2}). In order to overcome this difficulty, we choose another convolution approach.

Let the arithmetic function $u_k(a;c;n)$ be defined by
\begin{eqnarray}
u_k(a;c;n):=\sum_{n=l^ad^c}\mu_k(d),
\end{eqnarray}
which satisfies
\begin{eqnarray}
\sum_{n=1}^{\infty}\frac{u_k(a;c;n)}{n^s}=\frac
{\zeta(as)}{\zeta^k(cs)},\Re s>1.
\end{eqnarray}
We note that when $a=k=1,$ the function $u_k(a;c;n)$ is just the
characteristic function of the set of the $c$-free integers. Hence
we can write
\begin{equation}
f_k(a,b;c;n)=\sum_{n=m_1m_2^b}u_k(a;c;m_1).
\end{equation}

The function $u_k(a;c;n)$ plays an important role in this case.
Define
\begin{eqnarray}
\Delta_k(a;c;x):=\sum_{n\le x}u_k(a;c;n)-\frac {x^{\frac
1a}}{\zeta^k(\frac ca)}:=A_k(a;c;x)-\frac {x^{\frac
1a}}{\zeta^k(\frac ca)}.
\end{eqnarray}

For $\Delta_k(a;c;x)$, similar to Theorem 1 we also have

\begin{thm}
Let $x\ge 2$, $a<c$ and $\Delta_k(a;c;x)$ be defined by (1.17). If
the RH is true, then for any $ 1\le y< x^{\frac 1c}$ we have
\begin{eqnarray}
\Delta_k(a;c;x)=\sum_{l\le y}\mu_k (l)\psi\left((\frac
{x}{l^c})^{\frac 1a}\right)+O\left(x^{\frac
{1}{2a}+\varepsilon}y^{\frac 12-\frac {c}{2a}}+y^{\frac
12+\varepsilon}\right).
\end{eqnarray}
\end{thm}

On taking $y=x^{\frac {1}{a+c}}$ in Theorem 2 we get immediately
the following

\begin{cor}  Under the conditions of Theorem 2, we have
\begin{eqnarray}
\Delta_k(a;c;x)\ll x^{\frac {1}{a+c}+\varepsilon}.
\end{eqnarray}
\end{cor}

Now we state our main result for the case $c<b$, which improves
Theorem 1 in the case $a<c<b<2c$.

\begin{thm} Suppose RH is true.
 Let $x\ge 2$, $a<c<b<2c$,  $\Delta_k(a,b;c;x)$ and $\Delta_k(a,c;x)$ be defined
by (1.3) and (1.17) respectively. Suppose  $\Delta_k(a;c;x)\ll
x^{\alpha_k{(a;c)}+\varepsilon}$ such that $\alpha_k{(a;c)}<1/b$(a
natural restriction). Then for any $1\le y<x^{\frac 1b}$ we have
\begin{align}
\Delta_k(a,b;c;x)&=\sum_{d\le y}\Delta_k\left(a;c;\frac
{x}{d^b}\right)- \sum_{m\le \frac
{x}{y^b}}u_k(a;c;m)\psi\left((\frac x m)^{\frac
1b}\right)\\
&\ \ \  +O\left(x^{\frac {1}{2c}}y^{1-\frac
{b}{2c}}+(xy^{-b})^{\alpha_k{(a;c)}}+x^{\frac 1a}y^{-1-\frac
ba}\right).\nonumber
\end{align}
\end{thm}

\begin{cor}
Under the conditions of Theorem 3, we have
\begin{eqnarray}
\Delta_k(a,b;c;x)\ll x^{\frac {1}{a+b-ab\alpha_k
(a;c)}+\varepsilon}.
\end{eqnarray}
\end{cor}

{\bf Remark 1.3.} Since we use different convolution approaches in
Theorem 1 ($b<c$) and Theorem 3($b>c$), the exponential sums appeared   in these two theorems are also different. Hence we have to use different ways to estimate exponential sums in these two theorems. We note that  Corollary 3 implies
$\frac {1}{a+b-ab\alpha_k
(a;c)} <\frac 1b$, hence in the asymptotic formula (1.3), the second main term becomes a real main term.

{\bf Remark 1.4.} All corollaries above can be further improved by more  precise estimate for the exponential sums
involved(e.g., see \cite{ba2,bp,cz1,fi,hu1,mn,rs,sw1,wu2}).

 \bigskip

 The organization of this paper is as follows. In section 2 we shall give short proofs of  Theorem 1 and Theorem 2. The proof of Theorem 3 will be given in section 3.   In section 4, by the well-known Heath-Brown's method  we shall further improve Corollary 2 and obtain  a non-trivial estimate for $\Delta_k(a;c;x) (k=1,2)$, and then give some of its  applications to problems related to  the exponential convolution.
 In section 5 we discuss some applications of Corollary 1.1 and Corollary 1.3.
Finally, in section 6 we give an example to explain how to  get a
sharper upper bound by Theorem 3.

\bigskip

{\bf Notation.}  Throughout this paper $ \varepsilon$ denotes a
fixed positive constant, not necessarily the same in all
occurrences. As usual, let $\tau(n)$ and $\omega(n)$ denote the
divisor function , and the number of prime factors of $n$,
respectively. We also use $\tau_k(n)$ to denote the number of
decompositions of $n$ into $k$ factors, and let $\tau_1(n)=1$. Let
$q_r(n)$ denote the characteristic function of the set of $r$-free
integers. $x> 1$ is real, $\mathcal{L}=\log x$, $\{t\}$ denotes
the fractional part of $t, \psi(t)=\{t\}-1/2, \Vert t\Vert=\min
(\{t\},1-\{t\})$. We let $e(t)=\exp (2\pi it)$ and  $\delta
(x)=\exp(-A (\log x)^{\frac 35}(\log \log x)^{-\frac 15})$ for
some fixed constant $A>0$. $m\sim M$ means that $cM<m<CM$ for some
constants $0<c<C$.

\section{The proof of Theorem 1 and 2}

\begin{lem}
Let
 $A(s)=\sum_{n=1}^{\infty}a(n)n^{-s}$ converge absolutely for $\Re s=\sigma>\sigma_a$, and let
functions $H(u)$ and $B(u)$ be monotonically increasing such that
$$|a(n)|\le H(n),\ \  (n\ge
 1),$$
 $$\sum_{n=1}^{\infty}|a(n)|n^{-\sigma}\le B(\sigma),\ \ \sigma>\sigma_a.$$
If $s_0=\sigma_0+it_0, b_0>\sigma_a, b_0\ge b>0,
b_0\ge\sigma_0+b>\sigma_a,  T\ge
 1$ and $x\ge 1$, then for $x\notin \mathbb{N}$
\begin{align*}
\sum_{n\le x}a(n)n^{-s_0}=&\frac {1}{2\pi
i}\int_{b-iT}^{b+iT}A(s_0+s)\frac{x^s}{s}\mathrm{d}s+O\left(\frac{x^bB(b+\sigma_0)}{T}\right)\\
&+O\left(x^{1-\sigma_0}H(2x)\min (1,\frac {\log
x}{T})\right)+O\left( x^{-\sigma_0}H(N)\min (1,\frac {x}{T\Vert
x\Vert})\right).
\end{align*}
\end{lem}
\proof This lemma is the  well-known Perron's formula, for example,
see Theorem 2 of  page 98 in Pan\cite{pa}.\qed

Let $y\ge 1$ and define
\begin{eqnarray}
g_y(s):=\sum_{n>y}\frac {\mu_k(n)}{n^s}, \ \ \ \Re s>1.
\end{eqnarray}

\begin{lem} Suppose RH is true, then
$g_y(s)$ can be continued analytically to $\Re s=\sigma>\frac
12+\varepsilon$, and we have uniformly for $\sigma$ that
\begin{eqnarray}
 g_y(s)\ll y^{\frac 12 -\sigma + \varepsilon}(|t| +1)^{\varepsilon}, \sigma\ge
 \frac 12+ \varepsilon.
\end{eqnarray}
\end{lem}
\proof This lemma follows from Lemma 3 of Nowak\cite{nowak} immediately.
.\qed

{\bf The Proof of Theorem 1 and 2.}
 Let $\delta=\frac{\varepsilon}{10}$
  and $1\le y< x^{\frac 1c}$. We begin the proof of Theorem 1 in the same way as that
of Theorem 1 in Montgomery and Vaughan\cite{mv}. Here we only give
the details of  our proof for the case $a<b$. The proof  for the
case $a=b$  is similar.

Define
\begin{eqnarray}
f_{1,y}(n):=\sum_{ \stackrel { l^cm=n}{ l\le y}
 } \mu_k(l)d(a,b;m),
 f_{2,y}(n):=\sum_{ \stackrel { l^cm=n}{ l> y}
 } \mu_k(l)d(a,b;m),
\end{eqnarray}
hence
\begin{eqnarray}
f_k(a,b;c;n)=\sum_{ l^cm=n}
\mu_k(l)d(a,b;m)=f_{1,y}(n)+f_{2,y}(n).
\end{eqnarray}
We now write $A_k(a,b;x)$ in the form
\begin{eqnarray}
A_k(a,b;c;x):=\sum_{n\le x}f_k(a,b;c;n)=S_1(x)+S_2(x),
\end{eqnarray}
where
\begin{eqnarray}
S_1(x)=\sum_{n\le x}f_{1,y}(n),
\end{eqnarray}
and
\begin{eqnarray}
S_2(x)=\sum_{n\le x}f_{2,y}(n).
\end{eqnarray}
We first evaluate $S_1(x)$. From (1.8) we get

\begin{align}
S_1(x)&= \sum_{ \stackrel { l^cm\le x}{ l\le y}
 } \mu_k(l)d(a,b;m)=\sum_{  l\le y
 } \mu_k(l)D\left(a,b;\frac {x}{l^c}\right)\\
&=\sum_{l\le y}\mu_k (l)\left(\frac {\zeta(\frac ba)}{l^{\frac
ca}}x^{\frac 1a}+\frac
{\zeta(\frac ab)}{l^{\frac cb}}x^{\frac 1b}+\Delta\left(a,b;\frac {x}{l^c}\right)\right)\nonumber\\
&=\zeta (\frac ba)x^{\frac 1a}\sum_{l\le y}\frac
{\mu_k(l)}{l^{\frac ca}}+\zeta (\frac ab)x^{\frac 1b}\sum_{l\le
y}\frac {\mu_k(l)}{l^{\frac cb}}+\sum_{l\le
y}\mu_k (l)\Delta\left(a,b;\frac {x}{l^c}\right).\nonumber\\
&=\Res \left( \zeta (as)\zeta (bs)x^ss^{-1}\sum_{l\le y}\frac
{\mu_k(l)}{l^{cs}},\frac 1a \right)\nonumber \\
&\ \ \ +\Res \left( \zeta (as)\zeta (bs)x^ss^{-1}\sum_{l\le
y}\frac {\mu_k(l)}{l^{cs}},\frac 1b \right) +\sum_{l\le y}\mu_k
(l)\Delta\left(a,b;\frac {x}{l^c}\right).\nonumber
\end{align}

From (2.1) we have
\begin{eqnarray}
\sum_{n=1}^{\infty}\frac {f_{2,y}(n)}{n^s}=g_y(cs)\zeta(as)\zeta
(bs).
\end{eqnarray}

From (2.3), (2.7), (2.9) and Lemma 2.1 we obtain that
\begin{eqnarray}
S_2(x)=\frac {1}{2\pi i}\int_{\frac 1a +\varepsilon-ix^2}^{\frac
1a+\varepsilon+ix^2}g_y(cs)\zeta (as)\zeta(bs)x^ss^{-1}\mathrm{d}s
+O(x^\delta),
\end{eqnarray}
since $f_{2,y}(n)\ll n^\delta$ by a divisor argument.

{\bf Case (i).} $a<b<2a$. When we move the line of integration to
$\Re s=\sigma =\frac {1}{2a}+\delta_0$ with $\delta_0=\min
\{\delta, \frac {1}{2b}(1-\frac {b}{2a})\}$, then by the residue
theorem
\begin{align}
&\frac {1}{2\pi i}\int_{\frac 1a +\varepsilon-ix^2}^{\frac
1a+\varepsilon-ix^2}g_y(cs)\zeta (as)\zeta(bs)x^ss^{-1}\mathrm{d}s\\
&=\Res\left( g_y(cs)\zeta (as)\zeta (bs)x^ss^{-1},\frac 1a
\right)+\Res\left( g_y(cs)\zeta (as)\zeta (bs)x^ss^{-1},\frac 1b
\right)\nonumber\\
&\ \ \ +I_1+I_2-I_3,\nonumber
\end{align}
where
\begin{align*}
&I_1=\frac {1}{2\pi i}\int_{\frac {1}{2a}+\delta_0+ix^2}^{\frac
1a+\varepsilon+ix^2}g_y(cs)\zeta (as)\zeta
(bs)x^ss^{-1}\mathrm{d}s, I_2=\frac {1}{2\pi i}\int_{\frac
{1}{2a}+\delta_0-ix^2}^{\frac {1}{2a}+\delta_0 +ix^2}g_y(cs)\zeta
(as)\zeta (bs)x^ss^{-1}\mathrm{d}s,\\
&I_3=\frac {1}{2\pi i}\int_{\frac {1}{2a}+\delta_0 -ix^2}^{\frac
1a+\varepsilon-ix^2}g_y(cs)\zeta (as)\zeta
(bs)x^ss^{-1}\mathrm{d}s.
\end{align*}

From Lemma 2.2, we have
$$g_y(cs)\ll y^{\frac 12-\frac {c}{2a}}(|t|^\delta +1),( \sigma \ge \frac {1}{2a} +\delta)$$
Thus
\begin{eqnarray}g_y(cs)\zeta (as)\zeta (bs)\ll
y^{\frac 12-\frac {c}{2a}}(|t|^{3\delta} +1),( \sigma \ge \frac
{1}{2a} +\delta).
\end{eqnarray}
From (2.12) it is not difficult to see that
\begin{eqnarray}
I_j\ll y^{\frac 12-\frac {c}{2a}}x^{\frac
{1}{2a}+8\delta},(j=1,2,3).
 \end{eqnarray}

Now combining (1.3), (2.5), (2.8), (2.10), (2.11) and (2.13)
completes the proof of Theorem 1 in this case.

{\bf Case (ii).}  $b\ge 2a$. In this case, moving the line of
integration in (2.10) to $\Re s=\sigma =\frac {1}{2a}+\delta$, we
can treat $S_2(x)$ as in the above case  except the second residue
in relation (2.11) vanishes. In addition, applying Abel summation
formula and the estimate $\sum_{n\le x}\mu_k(n)\ll x^{\frac 12
+\varepsilon}$(this can be proved in the same way as that of Theorem
14.25(C) in Titchmarsh\cite{ti}, also see (2.10) in Nowak\cite{nowak}), it is easy to check that
\begin{eqnarray*}
x^{\frac 1b}\sum_{l>y}\frac {\mu_k(l)}{l^{\frac cb}}\ll x^{\frac
1b+\delta}y^{\frac 12-\frac cb}\ll x^{\frac
{1}{2a}+\delta}y^{\frac {1}{2}-\frac {c}{2a}}.
\end{eqnarray*}
Hence
\begin{align}
x^{\frac 1b}\sum_{l\le y}\frac {\mu_k(l)}{l^{\frac cb}}&=x^{\frac
1b}\sum_{l=1}^{\infty}\frac {\mu_k(l)}{l^{\frac cb}}-x^{\frac
1b}\sum_{l> y}\frac {\mu_k(l)}{l^{\frac
cb}}\\
&=\frac {1}{\zeta^k(\frac cb)}x^{\frac 1b}+O\left( x^{\frac
{1}{2a}+\delta}y^{\frac {1}{2}-\frac {c}{2a}}\right)\nonumber.
\end{align}

Therefore, Theorem 1 also holds in this case.  This completes the
proof of  Theorem 1.

The proof of Theorem 2 is very similar to that of Theorem 1, we omit
the details here. \qed

\section{\bf The proof of Theorem 3}
\begin{lem}
Let $a<c<b$ and $\Delta_k(a;c;x)$ be defined by (1.17). If
$\Delta_k(a;c;x)\ll x^{\alpha{(a;c)}+\varepsilon}$ such that
$\alpha{(a;c)}<  1/b$, then for $s>\alpha(a;c)$ we have
\begin{align}
\sum_{m\le x}u_k(a;c;m)m^{-s}=&\frac {x^{\frac
1a-s}}{(1-as)\zeta^k(\frac ca)}+\frac {\zeta
(as)}{\zeta^k(cs)}\\
&+\Delta_k(a;c;x)x^{-s}-s\int_{x}^{\infty}\Delta_k(a;c;t)t^{-s-1}\mathrm{d}t.\nonumber
\end{align}
\end{lem}

\proof By partial summation formula and (1.17) we get
\begin{align}
&\sum_{m\le
x}u_k(a;c;m)m^{-s}\\
&=A_k(a;c;x)x^{-s}+s\int_{1}^{x}A_k(a;c;t)t^{-s-1}\mathrm{d}t\nonumber\\
&=\frac {x^{\frac 1a-s}}{\zeta^k(\frac
ca)}+\Delta_k(a;c;x)x^{-s}+s\int_{1}^{x}\left(\frac {t^{\frac
1a}}{\zeta^k(\frac
ca)}+\Delta_k(a;c;t)\right)t^{-s-1}\mathrm{d}t\nonumber\\
&=\frac {x^{\frac 1a-s}}{(1-as)\zeta^k(\frac
ca)}+\Delta_k(a;c;x)x^{-s}-\frac {as}{(1-as)\zeta^k(\frac
ca)}+s\int_{1}^{x}\Delta_k(a;c;t)t^{-s-1}\mathrm{d}t.\nonumber
\end{align}

Suppose that $s>1$, we have from (1.15) and the condition
$\alpha{(a;c)}<\frac 1b$, when $x\rightarrow \infty$
\begin{eqnarray}
\frac {\zeta (as)}{\zeta^k(cs)}=-\frac {as}{(1-as)\zeta^k(\frac
ca)}+s\int_{1}^{\infty}\Delta_k(a;c;t)t^{-s-1}\mathrm{d}t.
\end{eqnarray}
By analytic continuation this equation also holds for
$s>\alpha(a;c)$. Substituting (3.3) into (3.2)  completes the proof
of Lemma 3.1.\qed

\begin{lem} Let $x\ge 2$, $a<c<b$, and
$\Delta_k(a,b;c;x)$ be defined by (1.3). If $\Delta_k(a;c;x)\ll
x^{\alpha{(a;c)}+\varepsilon}$ such that $\alpha{(a;c)}<  1/b$, then
for any $1\le y<x^{  1/b}$ we have

\begin{align}
\Delta_k(a,b;c;x)=&\sum_{d\le y}\Delta_k\left(a;c;\frac
{x}{d^b}\right)- \sum_{m\le \frac
{x}{y^b}}u_k(a;c;m)\psi\left((\frac x m)^{\frac
1b}\right)\nonumber\\
& -\frac {x^{\frac 1b}}{b} \int_{\frac {x}{y^b}}^{\infty}\frac
{\Delta_k(a;c;t)}{t^{1+\frac 1b }}\mathrm{d}t+\psi
(y)\Delta_k\left(a;c; \frac {x}{y^b}\right)+O(x^{\frac
1a}y^{-1-\frac ba}).\nonumber
\end{align}
\end{lem}

\proof  Let $1\le y\le x^{1/b}$. Applying (1.1),(1.15)-(1.17) and
Dirichlet's hyperbolic argument, we get
\begin{align}
&A_k(a,b;c;x)=\sum_{n\le x}f_k(a,b;c;n)=\sum_{md^b\le x}u_k(a;c;m)\\
&=\sum_{d\le y}\sum_{m\le \frac {x}{d^b}}u_k(a;c;m)+\sum_{m\le
\frac
{x}{y^b}}u_k(a;c;m)\sum_{y<d\le (\frac xm)^{\frac 1b}}1\nonumber\\
&=\frac {x^{\frac 1 a}}{\zeta^k(\frac ca)}\sum_{d\le y}\frac
{1}{d^{\frac ba}}+x^{\frac 1b}\sum_{m\le \frac {x}{y^b}}\frac
{u_k(a;c;m)}{m^{\frac 1b}}-y\sum_{m\le \frac
{x}{y^b}}u_k(a;c;m)\nonumber\\
&\quad+\sum_{d\le y}\Delta_k\left(a;c;\frac {x}{d^b}\right)-
\sum_{m\le \frac {x}{y^b}}u_k(a;c;m)\psi\left((\frac x m)^{\frac
1b}\right)+ \psi (y)\sum_{m\le \frac {x}{y^b}}u_k(a;c;m).\nonumber
\end{align}

Applying Lemma 3.1 with $s=\frac 1b$, we get
\begin{align}
\sum_{m\le \frac {x}{y^b}}\frac {u_k(a;c;m)}{m^{\frac 1b}} =&\frac
{\left( \frac {x}{y^b}\right)^{\frac 1a-\frac 1b}}{(1-\frac
ab)\zeta^k(\frac ca)}+\frac {\zeta (\frac ab)}{\zeta^k(\frac
ca)}\\
&+\Delta_k(a;c; \frac {x}{y^b})\left( \frac
{x}{y^b}\right)^{-\frac 1b}-\frac 1b \int_{\frac
{x}{y^b}}^{\infty}\frac {\Delta_k(a;c;t)}{t^{1+\frac 1b
}}\mathrm{d}t.\nonumber
\end{align}

In addition, we have from the Euler-Maclaurin formula that
\begin{eqnarray}
\sum_{d\le y}\frac {1}{d^{\frac ba}}=\zeta (\frac ba)+\frac
{y^{1-\frac ba}}{(1-\frac ba)}-\psi(y)y^{-\frac ba}+O(y^{-1-\frac
ba}).
\end{eqnarray}

Applying (1.17) again we also have
 \begin{eqnarray}
 \sum_{m\le \frac
{x}{y^b}}u_k(a;c;m)=\frac {\left( \frac {x}{y^b}\right)^{\frac
1a}}{\zeta^k(\frac ca)}+\Delta_k(a;c; \frac {x}{y^b}).
\end{eqnarray}
Substituting (3.5)-(3.7) into (3.4), we get
\begin{align}
A_k(a,b;c;x)=&\frac {\zeta(\frac ba)}{\zeta^k(\frac ca)}x^{\frac 1
a}+\frac {\zeta(\frac ab)}{\zeta^k(\frac cb)}x^{\frac 1b}-
\sum_{m\le \frac {x}{y^b}}u_k(a;c;m)\psi\left((\frac x m)^{\frac
1b}\right)\\
&+\sum_{d\le y}\Delta_k\left(a;c;\frac {x}{d^b}\right) -\frac
{x^{\frac 1b}}{b} \int_{\frac {x}{y^b}}^{\infty}\frac
{\Delta_k(a;c;t)}{t^{1+\frac 1b }}\mathrm{d}t\nonumber\\
&+\psi (y)\Delta_k\left(a;c; \frac {x}{y^b}\right)+O(x^{\frac
1a}y^{-1-\frac ba}).\nonumber
\end{align}
Now Lemma 3.2 follows from (1.3) and (3.8) at once.\qed

\begin{lem} Let $\Delta_k(a;c;x) $ be defined by (1.17).
If  RH is true, then for any fixed $\delta>0$ we have
\begin{eqnarray}
\int_{1}^{T}\Delta_k(a;c;u)\mathrm{d}u\ll T^{1+\frac {1}{2c}+
\delta }.
\end{eqnarray}
\end{lem}

\proof It suffices   to prove that for any  $M>2$, we have
\begin{eqnarray}
\int_{M}^{2M}\Delta_k(a;c;u)\mathrm{d}u\ll M^{1+\frac {1}{2c}+
\delta}.
\end{eqnarray}

 Taking in Lemma 2.1  $H(n)=n^\varepsilon, B(\sigma)=(\sigma-1)^{-k}, b=1+1/\log M, T=M^5$, we
get
\begin{eqnarray*}
\sum_{n\leq x}u_k(a;c;u)=\frac {1}{2\pi i}\int_{b-iT}^{b+iT}\frac
{\zeta (as)}{\zeta^k(cs)}\frac{u^s}{s}\mathrm{d}s+O(M^{-4}).
\end{eqnarray*}

It is well-known that  if RH is true, then for any fixed
$0<\eta<1/2,$ we have
\begin{eqnarray}
 \zeta (s)\ll (|t|+1)^{\eta},\
\zeta^{-1}(s)\ll (|t| +1)^\eta,\ \ \sigma
>\frac 12 + \eta.
\end{eqnarray}
Moving the line of integration to $\Re s=\frac {1}{2c}+\delta$, we
have by (3.11) and the estimate $\zeta(s)\ll (1+|t|)^{1/2}\
(\sigma\geq 0)$ that
\begin{eqnarray*}
\Delta_k(a;c;u)=\frac {1}{2\pi i}\int_{\frac
{1}{2c}+\delta-iT}^{\frac {1}{2c}+\delta+iT}\frac {\zeta
(as)}{\zeta^k(cs)}\frac{u^s}{s}\mathrm{d}s+O(M^{-1}).
\end{eqnarray*}

Thus we have
\begin{align*}
\int_{M}^{2M}\Delta_k(a;c;u)\mathrm{d}u & = \int_{\frac
{1}{2c}+\delta-iT}^{\frac {1}{2c}+\delta+iT}\frac {\zeta
(as)}{\zeta^k(cs)}\frac {\mathrm{d}s}{s}\int_{M}^{2M}u^s \mathrm{d}u+O(1)\\
&=\int_{\frac {1}{2c}+\delta-iT}^{\frac
{1}{2c}+\delta+iT}\frac{\zeta
(as)\left((2M)^{1+s}-M^{1+s}\right)}{\zeta^k(cs)s(1+s)}\mathrm{d}s+O(1).\nonumber\\
&\ll M^{1+\frac
{1}{2c}+\delta}\int_{-T}^{T}\left|\frac{\zeta(a(\frac
{1}{2c}+\delta+it))}{\zeta^k(c(\frac
{1}{2c}+\delta+it))}\right|\frac{\mathrm{d}t}{(1+|t|)^2}+1\nonumber\\
&\ll M^{1+\frac {1}{2c}+  \delta   }.\nonumber
\end{align*}
Namely (3.10) holds. This completes the proof of Lemma 3.3.\qed

{\bf The Proof of Theorem 3.} Theorem 3 follows immediately from
Lemma 3.2 and Lemma 3.3 with $\delta =\frac
{\varepsilon}{10}$.\qed

{\bf The Proof of Corollary 1.3.} It is easy to check
\begin{eqnarray}
\sum_{d\le y}\Delta_k\left(a;c;\frac {x}{d^b}\right)\ll \sum_{d\le
y}\left(\frac {x}{d^b}\right)^{\alpha_k(a;c)}\ll
x^{\alpha_k(a;c)}y^{1-b\alpha_k(a;c)}
\end{eqnarray}
and
\begin{align}
\sum_{m\le \frac {x}{y^b}}u_k(a;c;m)\psi\left((\frac x m)^{\frac
1b}\right)&\ll \sum_{m\le \frac {x}{y^b}}|u_k(a;c;m)|\\
&\ll \sum_{m\le \frac {x}{y^b}}d(a,c;m)\ll \left(\frac
{x}{y^b}\right)^{\frac 1 a}.\nonumber
\end{align}
Taking $y=x^{\frac {1-a\alpha_k (a;c)}{a+b-ab\alpha_k (a;c)}}$ we
find that Corollary 1.3 is an immediate consequence of Theorem 3 and
the above two estimates.\qed

\section{\bf Estimates for $\Delta_k (a;c;x) (k=1,2)$ and an application of Theorem 2  }

\subsection{\bf Some preliminary lemmas}

To treat the exponential sums appeared in Theorem 1 and 2, for the
arithmetic function $\mu_k (k\ge 2)$ one needs an analogue of the well-known
Vaughan's identity of M\"obius function $\mu$. First we shall
prove such an identity.

\begin{lem} (Vaughan's identity). Let $1\le N_1<N$.
Suppose that $U, V $ be two parameters with $1\leq U, V\le N_1$.
Then for any arithmetic function $f$ we have
\begin{eqnarray}
\sum_{N_1<n\le N}\mu_k(n)f(n)={\sum}_1-{\sum}_2-{\sum}_3,
\end{eqnarray}
where
\begin{align*}
&{\sum}_1=\sum_{U<m\le N/V}A(m)\sum_{\stackrel{N_1/m<n\le N/m}{V<n}}\mu_k(n)f(mn),\\
&{\sum}_2=\sum_{U<m\le UV}B(m)\sum_{N_1/m<n\le N/m}\tau_k(n)f(mn),\\
&{\sum}_3=\sum_{m\le U}B(m)\sum_{N_1/m<n\le N/m}\tau_k(n)f(mn),\\
&A(m)=\sum_{\stackrel{ed_1=m}{e\le U}}\mu_k(e)\tau_{k}(d_1),\ \
B(m)=\sum_{\stackrel{d_1d_2=m}{d_1\le U,d_2\le
V}}\mu_k(d_1)\mu_k(d_2).
\end{align*}
\end{lem}

\proof  Let
\begin{eqnarray}
F(U,s):=\sum_{d\le U}\frac{\mu_k(d)}{d^s},\Re s>1,
\end{eqnarray}
then
\begin{align}
\frac {1}{\zeta^k(s)}&=\left(\frac
{1}{\zeta^k(s)}-F(V,s)\right)\left(1-\zeta^k(s)F(U,s)\right)-F(U,s)F(V,s)\zeta^k(s)\\
&\ \ +F(U,s)+F(V,s),\Re s>1.\nonumber
\end{align}
Equating coefficients from both sides of (4.3)  gives the following
identity
\begin{eqnarray}
\mu_k(n)=b_1(n)+b_2(n)+b_3(n)+b_4(n),
\end{eqnarray}
where
\begin{align*}
b_1(n)&=-\sum_{\stackrel{dm=n}{d>V,m>1}}\mu_k(d)\left(\sum_{\stackrel{ed_1=m}{e\le
U}}\mu_k(e)\tau_{k}(d_1)\right),\\
b_2(n)&=-\sum_{\stackrel{d_1d_2m=n}{d_1\le U,d_2\le
V}}\mu_k(d_1)\mu_k(d_2)\tau_{k}(m),\\
b_3(n)& =\begin{cases}
 \d  \mu_k(n)
      & \mbox{if $n\le U$,}\\[1em]
  \d  0 & \mbox{if $n>U,$}
 \end{cases}
\ \  b_4(n) =\begin{cases}
 \d  \mu_k(n)
      & \mbox{if $n\le V$,}\\[1em]
  \d  0 & \mbox{if $n>V.$}
 \end{cases}
\end{align*}

From (1.10) we have
\begin{equation}
\sum_{md=n}\tau_{k}(m)\mu_k(d) =\begin{cases}
 \d  1
      & \mbox{if $n=1$,}\\[1em]
  \d  0 & \mbox{if $n>1.$}
 \end{cases}
 \end{equation}
In the sum for $b_1(n)$ we can replace the condition $m>1$ by $m>U$,
since the sum over $m$ vanishes by (4.5) when $1<m\le U$. Now
multiplying the above identity (4.4)  by $f(n)$ we get (4.1).\qed

We will also exploit the following several lemmas. Lemma 4.2 is
Lemma 1 of Graham and Pintz\cite{gp}(also see Theorem 18 of
Vaaler\cite{va}), Lemma 4.3 is well-known, Lemma 4.4 is Lemma 6 of
Fouvry and Iwaniec\cite{fi}, Lemma 4.5 is Lemma 4 of
 the second paper in\cite{ca}(also see (2.1) in Wu\cite{wu2}), Lemma 4.6 is Lemma 12 of Cao\cite{ca}.

\begin{lem}  Suppose $H>0$. There is a function
$\psi^*(x)$ such that
\begin{align*}
\psi^*(x)&=\sum_{1\le |h|\le H}\gamma (h)e(hx), \ \ \gamma(h)\ll
\frac {1}{|h|}, \\
|\psi^*(x)-\psi (x)|&\le \frac {1}{2H+2}\sum_{|h|\le
H}\left(1-\frac {|h|}{H}\right)e(hx).
\end{align*}
\end{lem}

\begin{lem} Let $X\ne 0$ and $\nu\ne 0,1$. If
$(\kappa,\lambda)$ is an exponent pair, then
\begin{eqnarray*}
\sum_{n\sim N}e(Xn^{\nu})\ll \left(XN^{\nu-1}\right)^\kappa
N^\lambda+X^{-1}N^{-\nu+1}.
\end{eqnarray*}
\end{lem}

\begin{lem} Let $0<M\le N<\gamma N\le \lambda M$, and
$|a_n|\le 1$. Then we have
\begin{eqnarray*}
\sum_{N<n\le \gamma N}a_n=\frac
{1}{2\pi}\int_{-M}^{M}\left(\sum_{M<n\le \lambda
M}a_nn^{-it}\right)N^{it}(\gamma^{it}-1)t^{-1}\mathrm{d}t+O\left(\log
(2+M)\right).
\end{eqnarray*}
\end{lem}

\begin{lem} Let $x\ge 2$, $\alpha, \beta, \gamma$ be given
real numbers with $\alpha (\alpha -1)\beta\gamma\neq 0$,
$|a(m)|\le 1$, $b(n_1,n_2)\le 1$. Suppose $G=xM^\alpha N_1^\beta
N_2^\gamma$, $(\kappa, \lambda)$ is an exponent pair and
\begin{eqnarray*}
T(M,N_1,N_2)=\sum_{m\sim M}\sum_{n_1\sim N_1}\sum_{n_2\sim
N_2}a(m)b(n_1,n_2)e(xm^\alpha n_1^\beta n_2^\gamma).
\end{eqnarray*}
Then
\begin{eqnarray*}
T(M,N_1,N_2)\mathcal{L}^{-2}\ll
\left(G^{\kappa}M^{1+\lambda+\kappa}(N_1N_2)^{2+\kappa}\right)^{\frac
{1}{2+2\kappa}}+M^{\frac 12}N_1N_2+M(N_1N_2)^{\frac 12}+G^{-\frac
12}MN_1N_2.
\end{eqnarray*}
\end{lem}

\begin{lem} Let $x\ge 2$, $\beta, \gamma$ be given real
numbers with $\beta\gamma\neq 0$,$|a(m)|\le 1$, $|b(n)|\le 1$,
$(\kappa, \lambda)$ is an exponent pair. Suppose $D$ is a
subdomain of $\{(m,n):m\sim M, n\sim N\}$ bounded by finite
algebraic curves , $G=xM^\beta N^\gamma$ and
\begin{eqnarray*}
T_1(M,N)=\sum_{(m,n)\in D}a(m)b(n)\psi(xm^\beta n^\gamma).
\end{eqnarray*}
Then
\begin{eqnarray*}
T_1(M,N)\mathcal{L}^{-6}\ll
\left(G^{\kappa}M^{1+\lambda+\kappa}N^{2+\kappa}\right)^{\frac
{1}{2+2\kappa}}+M^{\frac 12}N+MN^{\frac 12}+G^{-\frac 12}MN.
\end{eqnarray*}
\end{lem}

\subsection{ \bf An estimate of $\Delta (a;c;x)$ }

Suppose $1\leq a<c$ are two fixed integers. In this subsection, we
shall estimate the error term $\Delta (a;c;x)$ defined by (1.17)
with $k=1.$

Let $r>1$ be a fixed real  number. The function $\Delta_k(1;r;x)$ is
defined on $[1,\infty)$  such that for any $1\leq y\leq x^{1/r}$ one
has
\begin{eqnarray}
\Delta_k(1;r;x) =\sum_{l\le y}\mu_k (l)\psi\left(\frac
{x}{l^r}\right)+O\left(x^{\frac {1}{2}+\varepsilon}y^{\frac 12-\frac
{r}{2}}+y^{\frac 12+\varepsilon}\right).
\end{eqnarray}

It follows easily from Theorem 2 that
\begin{eqnarray}
\Delta_k(a;c;x)=\Delta_k(1;\frac ca;x^{\frac 1a}).
\end{eqnarray}
Hence we only need to estimate $\Delta_k(1;r;x)$ for real $r>1$.

Now we define
\begin{align}
\alpha (r) =\begin{cases}
 \d  \frac{7}{8r+6}
      & \mbox{if $1< r\le 5$ and $r\ne 2$,}\\[1em]
\d  \frac{17}{54}
      & \mbox{if $r=2$,}\\[1em]
 \d  \frac{67}{514}      & \mbox{if $5<r\leq 6,$}\\[1em]
\d  \frac {11(r-4)}{12r^2-37r-41} & \mbox{if $6<r\le 12,$}\\[1em]
 \d  \frac {23(r-1)}{24r^2+13r-37} & \mbox{if $12< r\le 20.$}
 \end{cases}
\end{align}
For $r>20$, $\alpha (r)$ is defined by the following procedure.
Let $q\ge 2$ and $Q=2^q$. For every $r>20$, there is a unique
integer $q$ such that $\frac {12Q-q-5}{2}<r\le \frac
{24Q-q-6}{2}$. With this value of $q$, define
\begin{eqnarray}
\alpha (r)=\frac {(12Q-1)r-12Q+1}{12Qr^2+(6Qq+1)r-(6Qq+12Q+1)}.
\end{eqnarray}

In this subsection we shall prove that

\begin{thm} Let $\alpha (r)$ be defined  by (4.8) and (4.9),
respectively. If RH holds, then
$$\Delta(1;r;x)=O\left(x^{\alpha(r)+\varepsilon}\right).$$
\end{thm}

{\bf Remark 4.1.}  Certainly one can improve the exponent $\alpha
(r)$ further for some special values of $r$. For example, R. C.
Baker and K. Powell\cite{bp} obtained recently that  $\alpha
(3)=\frac {17}{74}$, $\alpha (4)=\frac {17}{94}$ and $\alpha
(5)=\frac {3}{20}$. In addition, for large values of $r$, one can
take $\alpha (r)=\frac {1}{r+c^*r^{1/3}}$ for some constant $c^*>0$
(see Theorem 2, \cite{gp}.)

From (4.7) and Theorem 4 we get

 \begin{cor}  Let $1\leq a<c$ be two
fixed integers. If RH holds, then
\begin{eqnarray}  \Delta(a;c;x)=O\left(x^{\frac {1}{a}\alpha(\frac
ca)+\varepsilon}\right).
\end{eqnarray}
\end{cor}

\proof    Theorem 4 is proved in  Jia\cite{ji} for the case $r=2$.
  S. W. Graham and J.
Pintz\cite{gp}   showed that Theorem 4 holds for any integer
$r>3$. However it is easily seen that the argument of \cite{gp}
can be applied to any $r\ge 2.$   So we only give a proof of Theorem
4 for $1<r<2$.

Taking $y=x^{\frac {4}{4r+3}}$ , by Theorem 2, (4.6) and a simple
splitting argument, an estimate for $x^{\alpha(r)+\varepsilon}\ll
Y\le y$
\begin{eqnarray}
\sum_{Y<l\le 2Y}\mu (l)\psi\left(\frac {x}{l^r}\right)\ll
x^{\alpha(r)+\varepsilon}
\end{eqnarray}
would suffice to complete the proof of Theorem 4.

Choose $U=Y^{\frac 12}, V=Y^{\frac 14}$. Let $|a(m)|\le 1$ and
$|b(n)|\le 1$ be any complex-valued arithmetic functions. If we can
show the estimates
\begin{eqnarray}
\sum_{U< m\le Y/V}a(m)\sum_{Y< mn \le 2Y}b(n)\psi\left(\frac
{x}{m^rn^r}\right)\ll x^{\alpha(r)+\varepsilon}
\end{eqnarray}
and
\begin{eqnarray}
\sum_{ m\le U}a(m)\sum_{Y< mn \le 2Y}\psi\left(\frac
{x}{m^rn^r}\right)\ll x^{\alpha(r)+\varepsilon},
\end{eqnarray}
then (4.11) follows from Lemma 4.1.

We first estimate the type II sum (4.12). Assume $N\ll M$, applying
Lemma 4.6 with $(\kappa,\lambda)=(\frac 12,\frac 12)$, we get
\begin{align}
&\mathcal{L}^{-6}\sum_{ m\sim M}a(m)\sum_{\stackrel {n\sim N}{Y<
mn \le
2Y}}b(n)\psi\left(\frac {x}{m^rn^r}\right)\\
&\ll \left(\left(\frac{x}{(MN)^r}\right)^{\frac 12}M^2N^{\frac
52}\right)^{\frac 13}+M^{\frac 12}N+MN^{\frac 12}+\left(\frac
{x}{(MN)^r}\right)^{-\frac 12}MN\nonumber\\
&\ll \left(x(MN)^{4-r}(MN)^{\frac 12}\right)^{\frac
{1}{6}}+(MN)N^{-\frac 12}+x^{-\frac 12}(MN)^{1+\frac r2}\nonumber\\
&\ll x^{\frac 16}Y^{\frac{9-2r}{12}}+YN^{-\frac 12}+x^{-\frac
12}Y^{1+\frac r2}\ll x^{\alpha(r)}.\nonumber
\end{align}
Hence (4.14) holds under the condition $N\ll M$. If $N\gg M$, using
Lemma 4.4 to separate the dependence between the variables $n$  and
$m$, then interchanging the roles of $m$ and $n$, we can show that
(4.14) also holds in this case. The estimate (4.12) follows from
(4.14) by a simple splitting argument.

Now we turn to estimate the type I sum (4.13). If $M\ge V$, by the
same  the argument as that of (4.14), we get
\begin{eqnarray}
\mathcal{L}^{-6}\sum_{ m\sim M}a(m)\sum_{\stackrel {n\sim N}{Y< mn
\le 2Y}}\psi\left(\frac {x}{m^rn^r}\right)\ll
x^{\alpha(r)}.
\end{eqnarray}

If $M\le V$, applying Lemma 4.2  with $H=Yx^{-\alpha (r)}$ and Lemma
4.3 with $(\kappa,\lambda)=(\frac 12,\frac 12)$, we get that
\begin{align}
&\sum_{ m\le V}a(m)\sum_{Y< mn \le 2Y}\psi\left(\frac
{x}{m^rn^r}\right)\\
& \ll \sum_{ m\le V}\left(\frac 1H\frac Ym+\sum_{1\le |h|\le
H}\frac {1}{|h|}\left|\sum_{Y/m< n \le 2Y/m} e\left(\frac {hx}{m^rn^r}\right)\right|\right)\nonumber \\
& \ll x^{\alpha (r)}\mathcal{L}+\sum_{ m\le V}\sum_{1\le h\le
H}\frac 1h\left( (hx)^{\frac 12}Y^{-\frac r2}+(hxm)^{-1}Y^{1+r}\right)\nonumber\\
& \ll x^{\alpha (r)}\mathcal{L}+x^{\frac 12}Y^{-\frac r2}H^{\frac
12}V+x^{-1}Y^{1+r}\mathcal{L}\nonumber\\
& \ll x^{\alpha (r)}\mathcal{L}+x^{\frac 12-\frac {\alpha (r)}{2}}Y^{\frac 34-\frac r2}+x^{-1}y^{1+r}\mathcal{L}\nonumber\\
&\ll x^{\alpha (r)}\mathcal{L}.\nonumber
\end{align}
(Here note that if $1.5\le r<2$, we use the bound $Y\gg
x^{\alpha(r)}$, otherwise we use $Y\le y$) Finally, it follows
from (4.15) and (4.16) that (4.13) always holds. This completes
the proof of Theorem 4.\qed

\subsection{\bf An estimate of $\Delta_2(a;c;x)$  }

 \begin{thm}
 Let $a,c$ be two fixed integers such that $1\leq a<c\leq 9a/2$  and
$\Delta_2(a;c;x)$ be defined by (1.17) with $k=2$. Assume  that RH
holds, then
\begin{eqnarray}
\Delta_2(a;c;x)\ll x^{ \frac{7}{8c+6a}+\varepsilon}.
\end{eqnarray}
\end{thm}

\proof  Similar to the proof of Theorem 4, we only need to show
that for $1< r \le  9/2$ one has
\begin{eqnarray}
\Delta_2(1;r;x)=O\left(x^{\beta(r)+\varepsilon}\right),\
\beta(r)=7/(8r+6).
\end{eqnarray}
Taking $y=x^{\frac {4}{4r+3}}$ , by Theorem 2, (4.6) and a simple
splitting argument, an estimate
\begin{eqnarray}
\sum_{Y<l\le 2Y}\mu_2 (l)\psi\left(\frac {x}{l^r}\right)\ll
x^{\beta(r)+\varepsilon}\ \ (x^{\beta(r)+\varepsilon}\ll Y\le y)
\end{eqnarray}
would suffice to complete the proof of Theorem 5.

Choose $U=Y^{\frac 12}, V=Y^{\frac 14}$.  Let $|a_1(m)|\le 1$ and
$|b_1(n)|\le 1$ be any complex-valued arithmetic functions. If we
can show
\begin{eqnarray}
\sum_{U< m\le Y/V}a_1(m)\sum_{Y< mn \le 2Y}b_1(n)\psi\left(\frac
{x}{m^rn^r}\right)\ll x^{\beta(r)+\varepsilon}
\end{eqnarray}
and
\begin{eqnarray}
\sum_{ m\le U}a_1(m)\sum_{Y< mn \le 2Y}\tau(n)\psi\left(\frac
{x}{m^rn^r}\right)\ll x^{\beta(r)+\varepsilon},
\end{eqnarray}
then (4.19) follows from Lemma 4.1.

The estimate (4.20) can be proved by the same approach of (4.12), so
we omit its detals. Hence we only need to prove (4.21). From (4.20)
we get easily that
$$\sum_{Y^{1/4}<m\le U}a_1(m)\sum_{Y< mn \le 2Y}\tau(n)\psi\left(\frac
{x}{m^rn^r}\right)\ll x^{\beta(r)+\varepsilon}. $$ So it suffices
for us  to prove
\begin{eqnarray}
\sum_{ m\le V}a_1(m)\sum_{Y< mn \le 2Y}\tau(n)\psi\left(\frac
{x}{m^rn^r}\right)\ll x^{\beta(r)+\varepsilon}.
\end{eqnarray}
Let $1\le M\le V$. Now we are in a position to estimate the
exponential sum
\begin{eqnarray}
S_r(M,Y):=\sum_{ M<m\le 2M}a_1(m)\sum_{Y< mn_1n_2 \le
2Y}\psi\left(\frac {x}{m^rn_1^rn_2^r}\right).
\end{eqnarray}

Without the loss  of generality, we suppose $n_1\ll n_2$, hence
$n_1\ll (YM^{-1})^{\frac 12}$. Applying a simple splitting
argument, we have for some $N_1\ll (YM^{-1})^{\frac 12}$
\begin{eqnarray}
\mathcal{L}^{-1}S_r(M,Y)\ll \sum_{M< m\le
2M}a_1(m)\sum_{N_1<n_1\le 2N_1}\sum_{Y< mn_1n_2 \le
2Y}\psi\left(\frac {x}{m^rn_1^rn_2^r}\right).
\end{eqnarray}
Now we consider two cases.

(Case i): $\frac 12\le N_1\ll VM^{-1}$. In this case,  applying
Lemma 4.2 with $H=Yx^{-\beta (r)}$ and Lemma 4.3 with
$(\kappa,\lambda)=(\frac 12,\frac 12)$, similar to the estimate of
(4.16), we can obtain
\begin{eqnarray}
S_r(M,Y) \ll x^{\beta(r)}\mathcal{L}^3.
\end{eqnarray}

(Case ii): $ VM^{-1}\le N_1\ll (YM^{-1})^{\frac 12}$. Applying Lemma
4.4 to separate the dependence between the variable $n_2$ and the
variables $m, n_1$, we get for $N_2=\frac {Y}{MN_1}$ that
\begin{align}
&\sum_{M< m\le 2M}a_1(m)\sum_{N_1<n_1\le 2N_1}\sum_{Y< mn_1n_2 \le
2Y}\psi\left(\frac {x}{m^rn_1^rn_2^r}\right)\\
&= \frac {1}{2\pi}\int_{-N_2}^{N_2}\left(\sum_{M< m\le
2M}a_1(m)\sum_{N_1<n_1\le 2N_1}\sum_{N_2<n_2\le
8N_2}n_2^{-it}\psi\left(\frac {x}{m^rn_1^rn_2^r}\right)\right)
 N_2^{it}(8^{it}-1)t^{-1}\mathrm{d}t+O\left(MN_1\mathcal{L}\right)\nonumber\\
&=\frac {1}{2\pi}\int_{-N_2}^{N_2}\left(\sum_{MN_1<d\le
4MN_1}\sum_{N_2<n_2\le 8N_2}c(d)n_2^{-it}\psi\left(\frac
{x}{d^rn_2^r}\right)\right)N_2^{it}\frac{8^{it}-1}{t}\mathrm{d}t+O\left(MN_1\mathcal{L}\right),\nonumber
\end{align}
where $$c(d)=\sum_{d=mn_1, M<m\le 2M, N_1<n_1\le 2N_1}a_1(m)\ll
d^\varepsilon.$$

If $N_2\gg MN_1$,   applying Lemma 4.6 with $(M,N)=(N_2,MN_1)$ and
$(\kappa,\lambda)=(\frac 12,\frac 12)$ to estimate the inner sum in
the above expression, we get(similar to (4.14))
\begin{align}
&\ \ \ \ x^{-\frac {\varepsilon}{2}}\sum_{MN_1<d\le
4MN_1}\sum_{N_2<n_2\le
8N_2}c(d)n_2^{-it}\psi\left(\frac {x}{d^rn_2^r}\right)\\
&\ll \left(\left(\frac {x}{Y^r}\right)^{\frac
12}N_2^2(MN_1)^{\frac 52}\right)^{\frac 13}+N_2^{\frac
12}(MN_1)+N_2(MN_1)^{\frac 12}+\left(\frac {x}{Y^r}\right)^{-\frac
12}MN_1N_2\nonumber\\
&\ll \left(x^{\frac 12}Y^{-\frac r2+2}(MN_1N_2)^{\frac
14}\right)^{\frac 13}+(MN_1N_2)^{\frac 12}(YM)^{\frac
14}+(MN_1N_2)(MN_1)^{-\frac 12}+x^{-\frac 12}Y^{\frac
r2+1}\nonumber\\
&\ll \left(x^{\frac 12}Y^{\frac {9-2r}{4}}\right)^{\frac
13}+Y^{\frac 12}(YV)^{\frac 14}+YV^{-\frac 12}+x^{-\frac
12}Y^{\frac r2+1}\ll x^{\beta(r)}.\nonumber
\end{align}
If $N_2\ll MN_1$, using the same approach but with $(M,N)=(MN_1,
N_2)$ in Lemma 4.6   we get that (4.27) still holds.

Combining   (4.26) and (4.27), we obtain that
\begin{eqnarray}
S_r(M,Y) \ll x^{\beta(r)+\varepsilon}
\end{eqnarray}
holds in the  Case ii.

The estimate (4.22) now follows from the proofs of the above two
cases. .\qed

\subsection{\bf An application of Theorem 2}

   The exponential convolution(e-convolution)  was introduced by M. V.
 Subbarao\cite{su1}. Let $n>1$ be an integer of canonical form $n=p^{a_1}_{1}\cdots
p^{a_s}_{s}$. An integer $d$ is called an exponential divisor
(e-divisor) of $n$ if $d=p^{b_1}_{1}\cdots p^{b_s}_{s}$, where
$b_1|a_{1},\cdots,b_s|a_{s}$. Let $\tau^{(e)}(n)$ denote the number
of exponential divisors of $n$, which is called the exponential
divisor function. Let $r\geq 2 $ be a fixed integer. The integer
$n>1$ is called exponentially $r$-free (e-$r$-free) if all the
exponents $a_1,\cdots, a_s$ are $r$-free. Let $q_r^{(e)}$ denote the
characteristic function of the set of e-r-free integers. The
e-unitary convolution was introduced by N. Minculete and L. T\'{o}th
\cite{mt}. The function $I(n)=1 (n\ge 1)$ has inverses
with respect to e-convolution and e-unitary convolution denoted by
$\mu^{(e)}(n)$ and $\mu^{(e)*}(n)$, respectively. These are the
unitary and exponential analogues of the M\"obius function. These
arithmetic functions attract the interests of many authors, see for
example
 \cite{cz2, ha, ks1, ks2, lu,mt, pw, sa1, sa2, sw2, ss,  to1, to2, to3}.

L. T\'{o}th \cite{to2} showed that the Dirichlet series of
$\mu^{(e)}$ is of the form
\begin{eqnarray}
\sum_{n=1}^{\infty}\frac{\mu^{(e)}(n)}{n^s}=\frac{\zeta
(s)}{\zeta^2(2s)}W_1(s), \Re s >1,
\end{eqnarray}
where $W_1(s):=\sum_{n=1}^{\infty}\frac{w_1(n)}{n^s}$ is
absolutely convergent for $\Re s>\frac {1}{5}$.

Let
\begin{eqnarray}
\Delta_{\mu^{(e)}} (x):=\sum_{n\le x}\mu^{(e)}(n)-A_1x,
\end{eqnarray}
where
\begin{eqnarray*}
A_1:=m(\mu^{(e)})=\frac{W_1(1)}{\zeta^2(2)}=\prod_{p}\left(1+\sum_{n=2}^{\infty}\frac
{\mu(n)-\mu(n-1)}{p^n}\right).\end{eqnarray*}

  L. T\'{o}th \cite{to2}
showed $ \Delta_{\mu^{(e)}} (x)=O\left(
x^{\frac{91}{202}+\varepsilon} \right)$ under RH. The exponent
$\frac{91}{202}$ was improved to $\frac {37}{94}$ by X. Cao and W.
Zhai\cite{cz2}.

Similarly,  N. Minculete and L. T\'{o}th\cite{mt} showed that the
Dirichlet series of $\mu^{(e)*}$ is of the form
\begin{eqnarray}
\sum_{n=1}^{\infty}\frac{\mu^{(e)*}(n)}{n^s}=\frac{\zeta
(s)}{\zeta^2(2s)}W_2(s), \Re s >1,
\end{eqnarray}
where $W_2(s):=\sum_{n=1}^{\infty}\frac{w_2(n)}{n^s}$ is
absolutely convergent for $\Re s>\frac {1}{4}$. Let
\begin{eqnarray}
\Delta_{\mu^{(e)*}} (x):=\sum_{n\le x}\mu^{(e)*}(n)-A_2x,
\end{eqnarray}
where
\begin{eqnarray*}
A_2:=m(\mu^{(e)*})=\frac{W_2(1)}{\zeta^2(2)}=\prod_{p}\left(1+\sum_{n=2}^{\infty}\frac
{(-1)^{\omega(n)}-(-1)^{\omega(n-1)}}{p^n}\right).\end{eqnarray*}

Under RH, the estimate  $ \Delta_{\mu^{(e)*}} (x)=O\left(
x^{\frac{91}{202}+\varepsilon} \right)$ was proved in \cite{mt}.

As an application of Theorem 2 and  Theorem 5,  from (4.29) and
(4.31) we get the following
\begin{thm}
 Let $ \Delta_{\mu^{(e)}} (x)$ and $\Delta_{\mu^{(e)*}} (x)$
be defined by (4.30) and (4.32), respectively. If   RH is true, then
we have
\begin{eqnarray}
\Delta_{\mu^{(e)}} (x)=O\left( x^{\frac{7}{22}+\varepsilon}
\right)
\end{eqnarray}
and
 \begin{eqnarray}
 \Delta_{\mu^{(e)*}} (x)=O\left(
x^{\frac{7}{22}+\varepsilon} \right).
\end{eqnarray}
\end{thm}

For comparison, we have numerically that
$$\frac{91}{202}=0.45049\cdots,\ \ \frac
 {37}{94}=0.39361\cdots,\ \ \frac
{7}{22}=0.31818\cdots.$$

\section{\bf Some applications of Corollary 1.1 and Corollary 1.3 }

\subsection{\bf The distribution of   generalized square-full integers}

\indent

In 1963, E. Cohen \cite{co1} generalized square-full integers in the
following way: Let $a $ and $b $ are fixed positive integers. Let
$n>1$ be an integer of canonical form $n=p^{a_1}_{1}\cdots
p^{a_r}_{r}$ and $R_{a,b}$ denote the set of all $n$ such that each
exponent $a_i (1\le i\le r)$ is either a multiple of $a$ or is
contained in the progression $at+b(t\ge 0).$ Obviously $R_{2,3}$ is
the set of square-full integers. Let $a\nmid b$ , $f_{(a,b)}(n)$
denote the characteristic function of the set $R_{a,b}$. By Lemma
2.1 in E. Cohen\cite{co1} one has
\begin{eqnarray}
\sum_{n=1}^{\infty}\frac{f_{(a,b)}(n)}{n^s}=\frac{\zeta (as)\zeta
(bs)}{\zeta(2bs)}, \Re s >1.
\end{eqnarray}
We are interested in the summatory function of $f_{a,b}(n).$

First consider  the case $a<b$. Suppose also that $a\nmid b$.  In
this case the problem is closely related to the estimate of
$\Delta(a,b;x).$ We take $a=a, b=b, c=2b$ and $k=1$ in (1.1), then
the estimate (1.4) implies $\Delta(a,b;2b;x)=\Omega(x^{\frac
{1}{2(a+b)}})$. Suppose $\Delta(a,b;x)\ll
x^{\alpha{(a,b)}+\varepsilon}$ such that $\alpha{(a,b)}<1/2b$. By
Corollary 1.1 with $(a,b,c,k)=(a,b,2b,1)$, under the RH we have the
asymptotic formula
\begin{eqnarray}
\sum_{n\le x}f_{(a,b)}(n)=\frac {\zeta (\frac ba)}{\zeta (\frac
{2b}{a})}x^{\frac 1a}+\frac {\zeta (\frac ab)}{\zeta ( 2)}x^{\frac
1b}+ O\left(x^{\frac {1-a\alpha (a,b)}{a+2b-4ab\alpha
(a,b)}+\varepsilon}\right),
\end{eqnarray}
which improves Theorem 3.2 of \cite{su2}.

The distribution of square-full integers (the case $a=2, b=3$) has
received special attention. In this special case, the error term
in (5.2) becomes $(x^{11/72+\varepsilon})$, which was first proved
in \cite{ns}. The exponent $11/72$ was improved by several
authors. The best known result is duo to Wu\cite{wu2}, who
obtained exponent $\frac{12}{85}=0.1411\cdots$(also subject to the
RH).

 Now we suppose $b<a.$  From (1.4) we get  $\Delta(b,a;2b;x)=\Omega(x^{\frac
{1}{4b}})$. So without the loss of generality, we always suppose
$b<a<4b.$

  If $b<a< 2b$.
Applying Corollary 1.1 with $(a,b,c,k)=(b,a,2b,1)$, under the RH
one has
\begin{eqnarray}
\sum_{n\le x}f_{(a,b)}(n)=\frac {\zeta (\frac ab)}{\zeta (
2)}x^{\frac 1b}+\frac {\zeta (\frac ba)}{\zeta (\frac
{2b}{a})}x^{\frac 1a}+O\left(x^{\frac {1-b\alpha
(a,b)}{3b-4b^2\alpha (a,b)}+\varepsilon}\right),
\end{eqnarray}
which improves Theorem 3.4 of \cite{su2}.

 Finally look at the case $2b<a<4b.$ In this case, D.
Suryanarayana\cite{su2} proved
 that(see Remark 3.3 therein)
 \begin{equation}
\sum_{n\le x}f_{(a,b)}(n)=\frac {\zeta (\frac ab)}{\zeta (
2)}x^{\frac 1b}+O(x^{1/2b} \delta(x))
 \end{equation}
unconditionally. If RH is true, then Remark 3.4 of D.
Suryanarayana\cite{su2} claimed that
\begin{eqnarray}
\sum_{n\le x}f_{(a,b)}(n)=\frac {\zeta (\frac ab)}{\zeta (
2)}x^{\frac 1b}+O\left(x^{\frac {2a-b}{5ab-4b^2}}\omega (x)\right),
\end{eqnarray}
where $\omega (x)=exp\{A\log x(\log\log x)^{-1}\}$, $A$ is a
positive absolute constant. Here we note that on the right-hand
side of (5.5) the second main term $\frac {\zeta (\frac ba)}{\zeta
(\frac {2b}{a})}x^{\frac 1a}$ is absorbed into the error term .

 Applying Corollary 1.3 with
$(a,b,c,k)=(b,a,2b,1)$ and $\Delta (b;2b;x)\ll x^{\frac
{17}{54b}+\varepsilon}$ in Corollary 4.1,  we have under RH that
\begin{eqnarray}
\sum_{n\le x}f_{(a,b)}(n)=\frac {\zeta (\frac ab)}{\zeta (
2)}x^{\frac 1b}+\frac {\zeta (\frac ba)}{\zeta (\frac
{2b}{a})}x^{\frac 1a}+O\left( x^{\frac
{54}{37a+54b}+\varepsilon}\right),
\end{eqnarray}
which took the second main term $\frac {\zeta (\frac ba)}{\zeta
(\frac {2b}{a})}x^{\frac 1a}$ out of the error term in (5.4) when
$2b<a<54b/17=3.176\cdots b.$

\subsection{\bf On the order of the error function of the $(l, r)$-integers}

 For given integers $l,r$ with $1<r<l$, we say an integer $n$ is a  $(l,
r)$-integers if it has  the form $m^ln$ where $m, n $ are integers
and $n$ is $r$-free.  The definition of the  $(l, r)$-integers was
introduced by M. V. Subbarao and V. C. Harris\cite{sh}. Let
$g_{l,r}(n)$ denote the characteristic function of the set of $(l,
r)$-integers. By Lemma 2.6 in M. V. Subbarao and D.
Suryanarayana\cite{ss1}  we  have
\begin{eqnarray}
\sum_{n=1}^{\infty}\frac{g_{l,r}(n)}{n^s}=\frac{\zeta (s)\zeta
(ls)}{\zeta(rs)}, \Re s >1.
\end{eqnarray}
Hence $g_{l,r}(n)=f_{1,l;r}(n).$  We define the error term by
\begin{eqnarray*}
\Delta(1,l;r;x):=\sum_{n\le x}g_{l,r}(n)-\frac {\zeta (l)}{\zeta
(r)}x-\frac {\zeta(\frac 1l)}{\zeta (\frac rl)}x^{\frac 1l}.
\end{eqnarray*}
From (1.4) we have $\Delta(1,l;r;x):=\Omega(x^{\frac {1}{2r}})$.

  If $l\ge 2r$, the distribution of $(l,r)$-integers is
almost the same as   the distribution of $r$-free numbers. From
Theorem 4 we get under RH that
\begin{eqnarray}
\sum_{n\le x}g_{l,r}(n)=\frac {\zeta (l)}{\zeta
(r)}x+O\left(x^{\alpha(r)+\varepsilon}\right),
\end{eqnarray}
where $\alpha(r)$ is defined by (4.8) and (4.9). In particular, if
$l\ge 4$ we have $\Delta(1,l;2;x)=O(x^{\frac
{17}{54}+\varepsilon})$.

  If $r<l< 2r$ and $\alpha(r)\geq 1/l$ we get that (5.8) holds too. However, if
$\alpha(r)<1/l,$ by Corollary 1.3 with $(a,b,c,k)=(1,l,r,1)$, we get
under RH that
\begin{eqnarray}
\sum_{n\le x}g_{l,r}(n)=\frac {\zeta (l)}{\zeta (r)}x+\frac
{\zeta(\frac 1l)}{\zeta (\frac rl)}x^{\frac 1l}+O\left(x^{\frac
{1}{l+1-l\alpha(r)}+\varepsilon}\right).
\end{eqnarray}
  In particular,
$\Delta(1,3;2;x)=O(x^{\frac {18}{55}+\varepsilon})$.

The previously best known error term is due to  M. V. Subbarao and
D. Suryanarayana\cite{ss2}. Since $\frac
{18}{55}=0.32727\cdots<\frac 13$, (5.8) and (5.9) answer a
conjecture of M. V. Subbarao and D. Suryanarayana in \cite{ss1}(see
page 123) for the special case $r=2$. It should be noted that if
$2\le r\le 10, l=r+1$, we get the second main term; if $r\ge 11,
r<l=r+2$, we also get the second main term. Hence (5.9) is an
substantial improvement to theirs.

\subsection{\bf The distribution of e-$r$-free integers}

In this subsection we consider the distribution of e-$r$-free
integers. For the distribution of e-square-free integers, J.
Wu\cite{wu1} showed that $\Delta_{q_2^{(e)}} (x)=O\left(x^{\frac
14}\delta (x)\right)$, improving an earlier result of M. V.
Subbarao\cite{su1}. In the general case,
 L. T\'{o}th\cite{to1} obtained  that $\Delta_{q_r^{(e)}} (x)=O\left(x^{\frac
{1}{2^r}}\delta(x)\right)$. Under RH,  X. Cao and W.
Zhai\cite{cz2} showed that $\Delta_{q_2^{(e)}} (x)=O(x^{\frac
{38}{193}+\varepsilon})$, improving the exponent $\frac 15$ of L.
T\'{o}th\cite{to2}. In this subsection we shall study this topic
more carefully.

\begin{thm} (i) The Dirichlet series of $q_r^{(e)}$ is of form
\begin{eqnarray}
Q_r^{(e)}(s):=\sum_{n=1}^{\infty}\frac{q_r^{(e)}(n)}{n^s}=\frac
{\zeta (s)\zeta((2^r+1)s)}{\zeta(2^rs)\zeta(2^{r+1}s)}U_r(s),\Re
s>1,
\end{eqnarray}
where $U_r(s)$ is absolutely convergent for $\Re s>\frac
{1}{2^{r+1}+1}$.

(ii) Let
\begin{eqnarray}
\Delta_{q_r^{(e)}} (x):=\sum_{n\le
x}q_r^{(e)}(n)-C_1(r)x-C_2(r)x^{\frac{1}{2^r+1}},
\end{eqnarray}
where
\begin{eqnarray}
 C_1(r):=\frac{\zeta(2^r+1)U_r(1)}{\zeta(2^r)\zeta(2^{r+1})},\ \
  C_2(r):=\frac{\zeta(\frac{1}{2^r+1})U_r(\frac{1}{2^r+1})}{\zeta(\frac{2^r}{2^r+1})\zeta(\frac{2^{r+1}}{2^r+1})}.
\end{eqnarray}
Then
\begin{eqnarray}
\Delta_{q_r^{(e)}} (x)=\Omega\left(x^{\frac {1}{2^{r+1}}}\right).
\end{eqnarray}

(iii) If the RH is true, we have
\begin{eqnarray}
\Delta_{q_r^{(e)}} (x)=O\left(x^{\frac
{1}{2^r+2-(2^r+1)\alpha(2^r)}+\varepsilon}\right),
\end{eqnarray}
where $\alpha(r)$ is defined by  (4.8) and (4.9). In particular,
$\Delta_{q_2^{(e)}} (x)=O(x^{\frac {38}{193}+\varepsilon})$, here
$\frac {38}{193}=0.1968\cdots<\frac 15$.
\end{thm}

{\bf Remark 5.1.} Note that we always have $\frac
{1}{2^r+2-(2^r+1)\alpha(2^r)}<\frac {1}{2^r+1}$. In addition when
$r=2$, if we use the new estimate $\alpha (4)=\frac {17}{94}$
proved by R. C. Baker and K. Powell\cite{bp} recently, one can
slightly improve the above result to $\Delta_{q_2^{(e)}}
(x)=O(x^{\frac {94}{479}+\varepsilon})$. In section 6 we shall
improve the exponent $\frac {94}{479}$ further by the exponential
sum method.

\proof  Since the function $q_r^{(e)}$ is multiplicative and
$q_r^{(e)}(p^\alpha)=q_r(\alpha )$ for every prime power
$p^\alpha$. For $r\ge 2$, it is easy to verify that
$q_r^{(e)}(p)=q_r^{(e)}(p^2)=\cdots=q_r^{(e)}(p^{2^r-1})=1$,
$q_r^{(e)}(p^{2^r})=0$,
$q_r^{(e)}(p^{2^r+1})=\cdots=q_r^{(e)}(p^{2^{r+1}-1})=1$, and
$q_r^{(e)}(p^{2^{r+1}})=0$. Hence for $\Re s>1$
\begin{eqnarray}
\sum_{n=1}^{\infty}\frac{q_r^{(e)}(n)}{n^s}=\prod_{p}\left(1+\sum_{m=1}^{\infty}\frac
{q_r(m)}{p^{ms}}\right)
\end{eqnarray}
Applying the product representation of Riemann zeta-function
\begin{eqnarray*}
\zeta
(s)=\prod_{p}(1+p^{-s}+p^{-2s}+p^{-3s}+\cdots)=\prod_{p}(1-p^{-s})^{-1},\Re
s>1,
\end{eqnarray*}
we have for $\Re s>1$
\begin{eqnarray}\zeta (s)\zeta
((2^r+1)s)=\prod_{p}\left((1-p^{-s})(1-p^{-(2^r+1)s})\right)^{-1}.
\end{eqnarray}
Let
\begin{align}
f_{q_r^{(e)}}(z):&=1+\sum_{m=1}^{\infty}q_r (m)z^m\\
&=1+z+\cdots+z^{2^r-1}+z^{2^r+1}+\cdots+z^{2^{r+1}-1}+
\sum_{m=2^{r+1}+1}^{\infty} q_r (m)z^m.\nonumber
\end{align}

 By a simple calculation one get for $|z|<1$
\begin{align*}
(1-z)(1-z^{2^r+1})&= 1 - z - z^{2^r+1} + z^{2^r+2},\\
f_{q_r^{(e)}}(z)(1-z)(1-z^{2^r+1})&=1-z^{2^r}-z^{2^{r+1}}+\sum_{m=2^{r+1}+1}^{\infty}c_mz^m,
\end{align*}
and
\begin{align*}
&\left(1+z^{2^r}+z^{2(2^r)}+z^{3(2^r)}+\cdots\right)\left(1+z^{2^{r+1}}+z^{2(2^{r+1})}+z^{3(2^{r+1})}+\cdots\right)\\
&=1+z^{2^r}+2z^{2(2^r)}+2z^{3(2^r)}+\cdots.
\end{align*}

 From the above two relations, we easily obtain for $ |z|<1$
\begin{eqnarray}
&&f_{q_r^{(e)}}(z)(1-z)(1-z^{2^r+1})(1-z^{2^r})^{-1}(1-z^{2^{r+1}})^{-1}\\
&&=1+\sum_{m=2^{r+1}+1}^{\infty}C_mz^m.\nonumber
\end{eqnarray}

Taking $z=p^{-s}$ in (5.18), then combining (5.15), (5.16) and
(5.17) completes the proof of (5.10) in Theorem 7.

Applying Theorem 2 of M. K\"{u}leitner and W. G. Nowak\cite{kn}
and (5.10), we immediately get (5.13). By Corollary 1.3 with
$(a,b,c,k)=(1,2^r+1,2^r,1)$ and Theorem 4, we obtain
\begin{eqnarray}
\Delta (1,2^r+1;2^r;x)=O\left(x^{\frac
{1}{2^r+2-(2^r+1)\alpha(2^r)}+\varepsilon}\right).
\end{eqnarray}
Now (5.14) follows from (5.10) , (5.19) and a simple convolution
argument at once, and this completes the proof of Theorem 7.\qed

\subsection{\bf The divisor problem over the set of $r$-free numbers}

 Let $r\geq 2$ be a fixed integer. Winfried  Recknagel\cite{re} and Hailiang Fen\cite{fen}
  investigated the divisor problem over the
set of $r$-free numbers. Hailiang  Fen\cite{fen} showed that
\begin{eqnarray}
\Psi_r(s):=\sum_{n=1}^{\infty}\frac{\tau (n)q_r(n)}{n^s}=\frac
{\zeta^2 (s)}{\zeta^{r+1}(rs)}V_r(s),\Re s>1,
\end{eqnarray}
where $V_r(s)$ is absolutely convergent for $\Re s>\frac
{1}{r+1}$. In particular, $V_2(s)=\zeta^2(3s)W_2(s)$, $W_2(s)$ is
absolutely convergent for $\Re s>\frac {1}{4}$.

Let
\begin{eqnarray}
\sum_{n\le x}\tau (n)q_r(n)=\Res_{s=1}\Psi_r(s)\frac
{x^s}{s}+\Delta_{(r)}^\tau(x).
\end{eqnarray}
Hence this problem is reduced to estimate the error term
$\Delta_{r+1}(1,1,r;x)$. Hailiang Fen\cite{fen} showed that
\begin{eqnarray}
\Delta_{r+1}(1,1,r;x)=\left\{\begin{array}{ll}
O(x^{1/r} \delta(x)),&\mbox{if $r=2,3$}\\
\Delta_{r+1}(1,1,r;x)=O(x^{\frac {131}{416}}(\log
x)^{\frac{26947}{8320}}) ,& \mbox{if $r\geq 4$},
\end{array}\right.
\end{eqnarray}where the second estimate in (5.22) follows from
M. N. Huxley's bound(see \cite{hu2})
\begin{equation}
\Delta (1,1;x)\ll x^{\frac{131}{416}}(\log
x)^{\frac{26947}{8320}}.\end{equation}

Applying Corollary 1.1 with $(a,b,c,k)=(1,1,2,3)$ and
$(a,b,c,k)=(1,1,3,4)$ respectively, and with the help of (5.23), we
obtain under RH that
\begin{eqnarray}
\Delta_{3}(1,1;2;x)=O(x^{\frac {285}{724}+\varepsilon}),\ \ \
\Delta_{4}(1,1;3;x)=O(x^{\frac {285}{878}+\varepsilon}).
\end{eqnarray}
Finally, from (5.20) and (5.24) we get immediately the following

\begin{thm} If RH is true, then
\begin{eqnarray}
\Delta_{(2)}^\tau(x)=O(x^{\frac {285}{724}+\varepsilon}),\ \ \
\Delta_{(3)}^\tau(x)=O(x^{\frac {285}{878}+\varepsilon}).
\end{eqnarray}
\end{thm}

\section{\bf The distribution of e-square-free integers}

In this section we shall use the method of exponential sums, and
$\alpha (4)=\frac {17}{94}$(recall Remark 4.1 in Section 4) proved
by R. C. Baker and K. Powell\cite{bp} to prove

\begin{thm}
 Let $\Delta_{q_2^{(e)}} (x)$ be defined by (5.11). Assume
that RH holds, then
\begin{eqnarray}
\Delta_{q_2^{(e)}} (x)=O\left(x^{\frac
{23}{124}+\varepsilon}\right).
\end{eqnarray}
\end{thm}

{\bf Remark 6.1.} For comparison, we have $\frac
{23}{124}=0.18548\cdots$ and $\frac {94}{479}=0.19624\cdots$.

In the proof of Theorem 9 we   need the following lemma (see Lemma
6.9 of Kr\"atzel\cite{kr}).

\begin{lem}
 Let $1\le Z<Z_1\le x^{\frac 1a}$, $a>0,b>0$.
If $(\kappa,\lambda)$ is any exponent pair and if
$$(2\lambda-1)a>2\kappa b,\ \ (2\lambda-1)b>2\kappa a,$$
then
\begin{eqnarray*}
\sum_{Z<n\le Z_1}\psi\left(\left(\frac {x}{n^a}\right)^{\frac
1b}\right)\ll x^{\frac {2(\kappa+\lambda-\frac 12)}{a+b}}\log
x+Z_1\left(\frac {Z_1^a}{x}\right)^{\frac 1b}\log x.
\end{eqnarray*}
\end{lem}

{\bf Proof of Theorem 9.} Similar to the proof of Theorem 7, we
need only to prove
\begin{eqnarray}
\Delta(1,5;4;x)=O\left(x^{\frac {23}{124}+\varepsilon}\right).
\end{eqnarray}
Let $\alpha =\frac {23}{124}$. Applying Theorem 3 with
$(a,b,c,k)=(1,5,4,1)$ and $y=x^{\frac {101}{744}}$, the following
two estimates
\begin{eqnarray}
\sum_{d\le y}\Delta (1,4;\frac {x}{d^5})\ll x^{\alpha+\varepsilon}
\end{eqnarray}
\begin{eqnarray}
\sum_{m\le \frac {x}{y^5}}u_1(1;4;m)\psi\left((\frac x m)^{\frac
15}\right)\ll x^{\alpha+\varepsilon}
\end{eqnarray}
would suffice to finish the proof of Theorem 9.

We first estimate the sum in (6.3). Let $y_1=x^{\frac {3}{62}}$,
we split the sum in (6.3) into two parts and write

\begin{eqnarray}
\sum_{d\le y}\Delta (1,4;\frac {x}{d^5})=\sum_{d\le y_1}\Delta
(1,4;\frac {x}{d^5})+\sum_{y_1<d\le y}\Delta (1,4;\frac
{x}{d^5}):=S_1+S_2.
\end{eqnarray}
Clearly, it follows from $\alpha (4)=\frac {17}{94}$ that
\begin{eqnarray}
S_1\ll \sum_{d\le y_1}\left(\frac {x}{d^5}\right)^{\frac
{17}{94}+\varepsilon}\ll x^{\frac {17}{94}+\varepsilon}y_1^{\frac
{9}{94}}\ll x^{\alpha+\varepsilon}.
\end{eqnarray}
 To estimate $S_2$, we discuss two cases.

{\bf Case (i)} $y_1\ll D\le y^{\frac {843}{10912}}=y_2$. Applying
Theorem 2 with $Y_*=x^{\frac 15}D^{-\frac {4}{5}}$ and a simple
splitting argument, we have for some $1\ll N\ll Y_*$
\begin{align}
\sum_{d\sim D}\Delta (1,4;\frac {x}{d^5})&=\sum_{d\sim
D}\left(\sum_{n\le Y_*}\mu (n)\psi\left(\frac
{x}{d^5n^4}\right)+O\left((\frac {x}{D^5})^{\frac
12+\varepsilon}{Y_*}^{-\frac
32}+{Y_*}^{\frac 12+\varepsilon}\right)\right)\\
&\ll \mathcal{L}\left|\sum_{d\sim D}\sum_{n\sim N}\mu
(n)\psi\left(\frac {x}{d^5n^4}\right)\right|+x^{\frac
15+\varepsilon}D^{-\frac {3}{10}}+x^{\frac
{1}{10}+\varepsilon}D^{\frac {3}{5}}\nonumber
\end{align}
Applying Lemma 4.6 with $(\kappa,\lambda)=(\frac 12,\frac 12)$ and
$(M,N)=(N,D)$, one get
\begin{align}
&\mathcal{L}^{-6}\sum_{d\sim D}\sum_{n\sim N}\mu
(n)\psi\left(\frac
{x}{d^5n^4}\right)\\
&\ll \left(\left(\frac {x}{D^5N^4}\right)^\kappa
N^{1+\kappa+\lambda}D^{2+\kappa}\right)^{\frac{1}{2+2\kappa}}+D^{\frac
12}N+DN^{\frac 12}+x^{-\frac 12}D^{\frac 72}N^{3}\nonumber\\
&\ll x^{\frac 16}+x^{\frac 15}D^{-\frac
{3}{10}}+x^{\frac {1}{10}}D^{\frac {3}{5}}+x^{\frac {1}{10}}D^{\frac {11}{10}}\nonumber\\
& \ll x^{\frac 15}y_1^{-\frac {3}{10}}+x^{\frac {1}{10}}y_2^{\frac
{11}{10}} \ll x^{\alpha+\varepsilon}.\nonumber
\end{align}
 Combining (6.7) and (6.8), in this case one obtain
\begin{align}
\sum_{d\sim D}\Delta (1,4;\frac {x}{d^5})\ll
x^{\alpha+\varepsilon}.
\end{align}
{\bf Case (ii)} $y_2=y^{\frac {843}{10912}}\ll D\le y$. Applying
Theorem 2 with $Y_*=x^{\frac {13}{62}}D^{-1}$, and a simple
splitting argument, we have for some $1\ll N\ll Y_*$
\begin{align}
\sum_{d\sim D}\Delta (1,4;\frac {x}{d^5})&=\sum_{d\sim
D}\left(\sum_{n\le Y_*}\mu (n)\psi\left(\frac
{x}{d^5n^4}\right)+O\left((\frac {x}{D^5})^{\frac
12+\varepsilon}{Y_*}^{-\frac
32}+{Y_*}^{\frac 12+\varepsilon}\right)\right)\\
&\ll \mathcal{L}\left|\sum_{d\sim D}\sum_{n\sim N}\mu
(n)\psi\left(\frac
{x}{d^5n^4}\right)\right|+x^{\alpha+\varepsilon}+x^{\frac
{13}{124}+\varepsilon}D^{\frac 12}.\nonumber
\end{align}
 Now applying lemma 4.2  with $H=DNx^{-\alpha }$ and
Lemma 4.3 with $(\kappa,\lambda)=BA^2BA(\frac 16,\frac 46)=(\frac
{13}{40},\frac {22}{40})$, we easily obtain that
\begin{align}
&\sum_{d\sim D}\sum_{n\sim N}\mu (n)\psi\left(\frac
{x}{d^5n^4}\right)\ll \sum_{n\sim N}\left| \sum_{d\sim
D}\psi\left(\frac {x}{d^5n^4}\right)\right|\\
&\ll \sum_{n\sim N}\left(\frac DH+\sum_{1\le |h|\le H}\frac
{1}{|h|}\left|\sum_{d\sim D}e\left(\frac
{hx}{d^5n^4}\right)\right|\right)\nonumber\\
&\ll x^{\alpha}+\sum_{n\sim N}\sum_{1\le h\le H}\frac 1h\left(
\left(\frac {hx}{n^4D^6}\right)^\kappa D^\lambda+\left(\frac
{hx}{n^4D^6}\right)^{-1}\right)\nonumber\\
&\ll x^{\alpha}+\left(\frac {x^{1-\alpha}}{N^3D^5}\right)^\kappa
D^\lambda N+x^{-1}N^5D^6\nonumber\\
&\ll x^{\alpha}+x^{\frac {13}{62}+(\frac
{23}{62}-\alpha)\kappa}D^{\lambda-2\kappa-1}+x^{\frac
{3}{62}}D\nonumber\\
&\ll x^{\alpha}+x^{\frac {1339}{4960}}y_2^{-\frac
{11}{10}}+x^{\frac {3}{62}}y \ll x^{\alpha}.\nonumber
\end{align}
(Here we use $1-3\kappa>0$.) Combining (6.10) and (6.11), one also
has
\begin{eqnarray} \sum_{d\sim D}\Delta (1,4;\frac {x}{d^5})\ll
x^{\alpha+\varepsilon}.
\end{eqnarray}

Combining the above two cases, then using a simple splitting
argument we obtain
\begin{eqnarray}
S_2 \ll x^{\alpha+\varepsilon}.
\end{eqnarray}

Hence (6.3) follows from (6.5), (6.6) and (6.13).

Now we turn to prove (6.4). From (1.14), applying the Drichlet's
hyperbolic argument, we have for any $1<Z< \frac {x}{y^5}=Y$
\begin{align}
&\sum_{m\le \frac {x}{y^5}}u_1(1;4;m)\psi\left((\frac x m)^{\frac
15}\right)=\sum_{m_1m_2^4\le Y}\mu
(m_2)\psi\left((\frac {x}{ m_1m_2^4})^{\frac 15}\right)\\
&=\sum_{m_1\le Z}\sum_{m_2^4\le \frac {Y}{m_1}}\mu
(m_2)\psi\left((\frac {x}{ m_1m_2^4})^{\frac 15}\right)
+\sum_{m_2\le (\frac{Y}{Z})^{\frac 14}}\mu (m_2)\sum_{Z<m_1\le
\frac {Y}{m_2^4}}\psi\left((\frac {x}{ m_1m_2^4})^{\frac
15}\right).\nonumber
\end{align}

Applying Lemma 6.1 with $(\kappa,\lambda)=A^2BA(\frac 16,\frac
46)=(\frac {2}{40},\frac {33}{40})$ and $(a,b,Z_1)=(1,5,\frac
{Y}{m_2^4})$, we have
\begin{eqnarray}
\sum_{Z<m_1\le \frac {Y}{m_2^4}}\psi\left((\frac {x}{
m_1m_2^4})^{\frac 15}\right)\ll \left(\frac
{x}{m_2^4}\right)^{\frac 18}\mathcal{L} +\frac
{Y}{m_2^4}\left(\frac {Y} {x}\right)^{\frac 15}\mathcal{L}.
\end{eqnarray}

On taking $Z=x^{\frac 17}Y^{-\frac 17}$, it follows from (6.14)
and (6.15)
\begin{align}
&\sum_{m\le \frac {x}{y^5}}u_1(1;4;m)\psi\left((\frac x m)^{\frac
15}\right)\\
&\ll \sum_{m_1\le Z}\left(\frac {Y}{m_1}\right)^{\frac
14}+\sum_{m_2\le (\frac{Y}{Z})^{\frac 14}}\left(\left(\frac
{x}{m_2^4}\right)^{\frac 18} +\frac {Y}{m_2^4}\left(\frac {Y}
{x}\right)^{\frac 15}\right)\mathcal{L}\nonumber\\
&\ll Y^{\frac 14}Z^{\frac 34}+x^{\frac 18}Y^{\frac 18}Z^{-\frac
18}\mathcal{L}+Y^{\frac 65}x^{-\frac 15}\mathcal{L}\nonumber\\
&\ll x^{\frac {3}{28}}Y^{\frac 17}\mathcal{L}+Y^{\frac
65}x^{-\frac 15}\mathcal{L}\ll x^{\frac 14}y^{-\frac
57}\mathcal{L}+xy^{-6}\mathcal{L}\ll x^{\alpha
}\mathcal{L}^2.\nonumber
\end{align}
This completes the proof of Theorem 9.\qed

\smallskip
\begin{flushleft}
Authors' addresses:

\smallskip

Xiaodong Cao\\
Dept. of Mathematics and Physics,\\
Beijing Institute of Petro-Chemical Technology,\\
Beijing, 102617, P. R. China\\
Email: caoxiaodong@bipt.edu.cn

\medskip

Wenguang Zhai,\\
Department of Mathematics, \\
China University of Mining and Technology, \\
Beijing 100083, P. R. China\\
E-mail:  zhaiwg@hotmail.com
\end{flushleft}

\end{document}